\documentclass[11pt]{amsart}
\usepackage{amsmath,amssymb,amsthm}
\usepackage[latin1]{inputenc}
\usepackage{version,tabularx,multicol}
\usepackage{graphicx,float}

\newcommand{\reff}[1]{(\ref{#1})}

\theoremstyle{plain}
\newtheorem{theo}{Theorem}[section]

\newtheorem{cor}[theo]{Corollary}
\newtheorem{prop}[theo]{Proposition}

\newtheorem{lem}[theo]{Lemma}

\theoremstyle{remark}
\newtheorem{rem}[theo]{Remark}

\newcommand{\cb}{{\mathcal B}}

\newcommand{\E}{{\mathbb E}}

\newcommand{\N}{{\mathbb N}}
\renewcommand{\P}{{\mathbb P}}

\newcommand{\R}{{\mathbb R}}

\newcommand{\ind}{{\bf 1}}

\newcommand{\val}[1]{\mathop{\left| #1 \right|}\nolimits}
\newcommand{\inv}[1]{\mathop{\frac{1}{ #1}}\nolimits}
\newcommand{\expp}[1]{\mathop {\mathrm{e}^{ #1}}}

\begin{document}

\title[CSBP with immigration]{Changing the branching
  mechanism of a continuous state branching process using  immigration}

\date{\today}
\author{Romain Abraham} 

\address{
Romain Abraham,
Université d'Orléans, Laboratoire MAPMO
CNRS, UMR 6628
Fédération Denis Poisson, FR 2964
Bâtiment de Mathématiques
BP 6759, 45067 Orleans cedex 2, France 
}
  
\email{romain.abraham@univ-orleans.fr}

\author{Jean-François Delmas}

\address{
Jean-Fran\c cois Delmas,
CERMICS, \'Ecole Nationale des Ponts et Chaussées, ParisTech, 6-8
av. Blaise Pascal, 
  Champs-sur-Marne, 77455 Marne La Vallée, France.}

\email{delmas@cermics.enpc.fr}

\begin{abstract}
  We consider  an initial population  whose size evolves according  to a
  continuous state  branching process.  Then we add  to this  process an
  immigration  (with  the  same   branching  mechanism  as  the  initial
  population), in such  a way that the immigration  rate is proportional
  to the whole population size. We prove this continuous state branching
  process  with immigration  proportional to  its own  size is  itself a
  continuous state branching process.  By considering the immigration as
  the apparition  of a new type,  this construction is a  natural way to
  model  neutral  mutation.  It  also  provides  in  some sense  a  dual
  construction of  the particular pruning  at nodes of  continuous state
  branching process introduced by the  authors in a previous paper.  For
  a  critical  or  sub-critical  quadratic branching  mechanism,  it  is
  possible  to  explicitly compute  some  quantities  of interest.   For
  example, we compute  the Laplace transform of the  size of the initial
  population conditionally on the non extinction of the whole population
  with  immigration.  We  also  derive the  probability of  simultaneous
  extinction  of the initial  population and  the whole  population with
  immigration.

\end{abstract}

\keywords{}

\subjclass[2000]{60G55, 60J25, 60J80, 60J85.}

\maketitle

\section{Introduction}

We consider an initial Eve-population of  type 0 whose size evolves as a
continuous state  branching process  (CB), $Y^0=(Y^0_t, t\geq  0)$, with
branching mechanism  $\psi^0$ defined by 
\begin{equation}
   \label{eq:def-psi0}
\psi^0 (\lambda)=\alpha^0\lambda+\beta \lambda^2+ \int_{(0,\infty)}\pi(d\ell)
\left[\expp{-\lambda\ell}-1+\lambda\ell\ind_{\{\ell \leq 1\}}\right],  
\end{equation}
where $\alpha^0\in  \R$, $\beta\geq 0$ and  $\pi$ is a  Radon measure on
$(0,\infty  )$  such  that  $\int_{(0,\infty  )} (1  \wedge  \ell^2)  \;
\pi(d\ell)<\infty $. See  \cite{l:csbp} for a definition of  CB as limit
of  Galton-Watson processes.  We assume  that this  population undergoes
some irreversible  mutations with constant rate, giving  birth either to
one individual  of type 1 (with  rate $\bar \alpha $),  or to infinitely
many  offsprings  of  type  1  (with rate  and  mutant  offsprings  size
described by a measure $\nu$).  This second population of type 1 evolves
according to  the same branching  mechanism as the  Eve-population (i.e.
the mutations  are neutral).  The  population of type 1  undergoes also
some mutations and  gives birth to a population of type  2 with the same
rules, and so on.

If we loose track of the genealogy,  the new population of type 1 can be
seen as an immigration process with rate proportional to the size of the
Eve-population, the population of type  2 is an immigration process with
rate proportional  to the size of the  population of type 1,  and so on.
We are interested in the law of the total population size $X=(X_t, t\geq
0)$, which is a CB with  immigration (CBI) rate proportional to its own size.
If the mutations are neutral, we expect  $X$ to be a CB.  This is indeed
the case : if $\psi^0$  is the branching mechanism of the Eve-population
and
\[
\phi(\lambda)=\bar \alpha
\lambda+\int_{(0,+\infty)}\nu(dx)\left(1-\expp{-\lambda 
  x}\right)
\]
is the  immigration mechanism,  then the total  population size is  a CB
with branching mechanism $\psi=\psi^0-\phi$, see Theorem \ref{th:main}.

Another  approach is  to associate  with critical  or sub-critical  CBs a
genealogical structure,  i.e. an infinite continuous  random tree (CRT),
see  \cite{lglj:bplpep}  or \cite{dlg:rtlpsbp}.  In  that context,  each
individual of the  CB can be followed during  its lifetime and mutations
can be added  as marks on its lineage.  Pruning  the CRT associated with
the  total  population (of  branching  mechanism $\psi=\psi^0-\phi$)  at
these  marks  allow  to   recover  the  Eve-population  from  the  total
population.  This  construction has been used in  \cite{ad:falp} and the
construction given here via immigration  proportional to the size of the
population  can be  seen  as the  dual  of the  pruning construction  of
\cite{ad:falp}, see Section  \ref{sec:dual} and more precisely Corollary
\ref{cor:dual}.   However \cite{ad:falp} considers  only the  case where
the branching  mechanism of  $X$ is  given by a  shift of  the branching
mechanism of the  Eve-population.  We shall give in  a forthcoming paper
\cite{adv:plcrt} a more general  pruning procedure which will correspond
to the general proportional immigration presented here.

Natural questions then  arise from a population genetics  point of view,
where only the whole population $X_t$  is observed at time $t$. In order
to  compute some  quantities related  to the  Eve-population,  given the
total population, we compute the joint  law of the Eve-population and the
whole population at a  given time: $(Y^0_t,X_t)$. For quadratic critical
or  sub-critical  branching   mechanism,  we  provide  explicit  Laplace
transform of  the joint distribution of $(Y^0_t,  X_t)$.  In particular,
we  compute $\P(Y^0_t=0| X_t>0)$,  the probability  for the  Eve-type to
have disappeared at time $t$, conditionally on the survival of the total
population at  time $t$, see Remark \ref{rem:Yt-cond}.   We also compute
the  Laplace transform  of $Y^0_t$  conditionally on  the  population to
never   be  extinct,   see  Proposition   \ref{prop:Yu_Xi}.    In  Lemma
\ref{lem:tau=sigma},   we  compute   the  probability   of  simultaneous
extinction  of the  Eve-population and  the whole  population,  in other
words, the probability  for the last individual alive  to have undergone
no  mutation. The techniques  used here  didn't lead  us to  an explicit
formula  but for  the quadratic  branching mechanism.   For  the general
critical  or sub-critical  case,  we use,  in \cite{ad:wdlcrtseppnm},  a
Williams' decomposition of  the genealogical tree to give  a very simple
formula  for   the  probability   of  simultaneous  extinction   of  the
Eve-population and the whole population.

In  the  particular  case  of  CB  with  quadratic  branching  mechanism
($\psi(u)=\beta   u^2$,  $\beta>0$),  similar   results  are   given  in
\cite{s:cddbt}   (using   genealogical   structure   for  CB)   and   in
\cite{w:bprtsbm}   (using  a  decomposition   of  Bessel   bridges  from
\cite{py:dbb}).   In  the   critical  ($\psi'(0^+)=0$)  or  sub-critical
($\psi'(0^+)>0$)  case  one could  have  used  the genealogical  process
associated to  CB introduced by \cite{lglj:bplpep} and  to CBI developed
by \cite{l:gcsbpi} to prove  the present result. This presentation would
have  been  more  natural  in   view  of  the  pruning  method  used  in
\cite{ad:falp}.   Our  choice  not  to  rely on  this  presentation  was
motivated   by  the   possibility  to   consider   super-critical  cases
($\psi'(0^+)<0$).

The  paper is  organized  as follows:  In  Section \ref{sec:CB-CBI},  we
recall   some  well   known   facts   on  CB   and   CBI.   In   Section
\ref{sec:(Y,fi)},  we  built a  CBI  $X$  whose  branching mechanism  is
$\psi^0$  and immigration  rate at  time $t$  proportional to  $X_t$ and
prove  that   this  process  is  again   a  CB.   We   give  in  Section
\ref{sec:dual} some links with the  pruning at nodes of CB introduced in
\cite{ad:falp}.    Eventually,  we   compute  the   joint  law   of  the
Eve-population and  the whole population in  Section \ref{sec:appli}, as
well as some related quantities.

\section{CB and CB with immigration}
\label{sec:CB-CBI}
The results from this section can be found in \cite{kw:bpirll} (see
also \cite{li:bpirt} for a survey on CB and CBI, and the references therein). 
Let $\psi $ be a branching mechanism of a CB: for $\lambda\geq 0$, 
\begin{equation}
   \label{eq:def_psi}
\psi (\lambda)=\alpha \lambda+\beta \lambda^2+ \int_{(0,\infty)}\pi(d\ell)
\left[\expp{-\lambda\ell}-1+\lambda\ell\ind_{\{\ell \leq 1\}}\right],  
\end{equation}
where $\alpha \in  \R$, $\beta\geq 0$ and  $\pi$ is a  Radon measure on
$(0,\infty  )$  such  that  $\int_{(0,\infty  )} (1  \wedge  \ell^2)  \;
\pi(d\ell)<\infty  $.  Notice $\psi $  is  smooth  on  $(0,\infty )$  and
convex.   We   have   ${\psi }'(0^+)\in   [-\infty   ,   +\infty   )$,   and
${\psi }'(0^+)=-\infty  $  if  and   only  if  $\int_{(1,\infty  )}  \ell  \;
\pi(d\ell)=\infty $. In order to
consider only conservative CB, we shall also assume
that for all $\varepsilon>0$
\begin{equation}
   \label{eq:conservatif}
\int_0^\varepsilon \inv{\val{\psi (u)}} \; du =\infty .
\end{equation}
Notice that $\psi'(0^+)>-\infty $ implies \reff{eq:conservatif}.


\subsection{CB}
\label{sec:CB}
Let $\P_x$ be the law of a  CB $Z=(Z_t,t\ge 0)$ started at $x\geq 0$ and
with  branching mechanism  $\psi $.  The process $Z$ is  a Feller process and thus càd-làg. Thanks to
\reff{eq:conservatif},   the   process    is   conservative,   that   is
a.s. for all $t\geq 0$, $Z_t<+\infty$.  
For every  $\lambda>0$, for
every $t\ge 0$, we have
\begin{equation}
\label{eq:laplace_csbp}
\E_x\left[\expp{-\lambda Z_t}\right]=\expp{-xu(t,\lambda)}
\end{equation}
where the function $u$ is the unique non-negative solution of
\begin{equation}
\label{eq:int_u}
u(t,\lambda)+\int_0^t\psi \bigl(u(s,\lambda)\bigr)ds=\lambda, \quad
\lambda\geq 0, \quad t\geq 0.
\end{equation}
This equation is equivalent to
\begin{equation}
   \label{eq:int_u2}
\int_{u(t,\lambda)}^\lambda \frac{dr}{\psi (r)}=t \quad
\lambda\geq 0, \quad t\geq 0.
\end{equation}

The  process $Z$  is infinitely  divisible.   Let $Q$  be its  canonical
measure.   The $\sigma$-finite  measure $Q$  is  defined on  the set  of
càd-làg  functions. Intuitively,  it gives  the ``distribution''  of the
size process for a  population generated by an infinitesimal individual.
In particular, under $\P_x$, $Z$  is distributed as $\sum_{i\in I} Z^i$,
where  $\sum_{i\in I}  \delta_{Z^i}$  is a  Poisson  point measure  with
intensity $x Q(dZ)$. Thus,  for any non-negative measurable function $F$
defined on the set of càd-làg functions, we have the following
exponential formula  for Poisson point measure
\[
\E_x[\expp{-\sum_{i\in I} F(Z^i)}]=\exp\left(- x Q[1-\expp{-F(Z)}]\right).
\]

The CB is called critical (resp. super-critical, resp. sub-critical) if
${\psi }'(0^+)=0$ (resp. ${\psi }'(0^+)<0$, resp. ${\psi }'(0^+)>0$). 

We shall  need inhomogeneous notation.  For $t<0$, we set  $Z_t=0$.  Let
$ \P_{x,t}$ denote the  law of $(Z_{s-t}, s\in \R)$  under $\P_x$, and let
$Q_t$ be the distribution of $(Z_{s-t}, s\in \R)$ under $Q$.

For $\mu$ a positive measure on $\R$, we set 
  $H^\mu=\sup\{r\in  \R;  \mu([r,\infty  ))>0\}$  the maximal element of
  its support.
\begin{prop}
\label{prop:uniq_x}
  Let $\mu$  be a finite positive  measure on $\R$  with  support
  bounded from above (i.e. $H^\mu$ is finite). Then we have for all
  $s\in \R$, $x\geq 0$,  
\begin{equation}
\label{eq:variante-mu}
\E_x\left[\expp{-\int Z_{r-s}\;  \mu(dr) }\right]=\expp{-xw(s)},
\end{equation}
where  the function  $w$ is  a measurable  locally  bounded non-negative
solution of the equation
\begin{equation}
 \label{eq:def_w}
w(s)+\int_s^{\infty }\psi (w(r))dr=\int_{[s,\infty )} \mu(dr), \quad s\leq
H^\mu\quad\text{and}\quad w(s)=0, \quad s> H^\mu. 
\end{equation}
If ${\psi }'(0^+)>-\infty   $ or if $\mu(\{H^\mu\})>0$, then \reff{eq:def_w}
has a unique measurable  locally  bounded non-negative
solution. 
\end{prop}
This result is well known for critical and sub-critical
branching mechanism (see eg. \cite{d:bpss}). As we didn't find a reference for
super-critical branching mechanism,  we give a short proof of
this Proposition. 

\begin{proof}
   
  Let        $n\geq        1$.        We       set        $\displaystyle
  Z^{(n),s}_t=Z_{\frac{i+1}{2^n}-s}$  for  $t\in  [i/2^n, (i+1)/  2^n)$.
  Using   that    $Z$   is   càd-làg,   we    get   a.s.   $\displaystyle
  \lim_{n\rightarrow  \infty }  Z_{t}^{(n),s}=Z_{t-s}$  for all  $t,s\in
  \R$. Since the process $Z$  is finite, we get by the dominated convergence
  theorem 
  a.s. for all $s\in \R$
\[
\int_{[-s, H^\mu]} Z_{r-s} \;\mu(dr)=\lim_{n\rightarrow \infty }
\int_{[-s, H^\mu]} Z_{r}^{(n),s} \;\mu(dr).
\]
Using the Markov property of $Z$, we get that 
\[
\E_x\left[\expp{-\int Z_{r}^{(n),s}\; \mu(dr)} \right]  =\expp{- x
    w^{(n)}(s)},
\]
where $w^{(n)}$ is the unique non-negative solution of
\[
w^{(n)}(s) + \int_s^{([H^\mu2^n]+1)/2^n} \!\!\!\psi (w^{(n)}(r)) \; dr=
\int_{[k/2^n, \infty )} \mu(dr), 
\]
with $k$ s.t. $k/2^n<s\leq  (k+1)/2^n$. 

Let  $T>  H^\mu  +1$.  Notice   that  for  all  $s\in  [-T,  T]$,  we
have $\int
Z_{r}^{(n),s}  \mu(dr)   \leq  \sup\{Z_t,  t\in  [0,2T]   \}  \mu([  -T,
H^\mu])<\infty   $   a.s.   Let    $C$   be   defined   by   $\expp{-C}=
\E_x[\expp{-\sup\{Z_t,  t\in  [0,2T]\} \mu([  -T,  H^\mu])} ]$.   Notice
$C<\infty $. This implies that for all $n\geq 1$, $ s\in [-T, T]$,
\[
0\leq w^{(n)}(s) \leq  C<\infty . 
\]
By   dominated   convergence    theorem,   $w^{(n)}(s)$   converges   to
$w(s)=-\log\left(\E_1[\expp{-\int Z_{r-s}  \mu(dr)}]\right)$, which lies
in $[0,C]$, for all $s\in  [-T,T]$. By dominated convergence theorem, we
deduce  that $w$ solves  \reff{eq:def_w}.  Since  $T$ is  arbitrary, the
Proposition  is   proved  but  for   the  uniqueness  of   solutions  of
\reff{eq:def_w}.

If  ${\psi }'(0^+)>-\infty   $,  then  $\psi $  is   locally Lipschitz. This
implies  there exists  a unique locally bounded non-negative solution  of
\reff{eq:def_w}. 

If  ${\psi }'(0^+)=-\infty $, and  $\mu(\{H^\mu\})>0$, we  get that  $ \int
Z_{r-s}  \;\mu(dr)\geq aZ_{H^\mu-s}$,  where  $a=\mu(\{H^\mu\})>0$. This
implies that $ w(s)\geq u(H^\mu -  s, a)>0$ for $s\in \R$.  The function
$u(\cdot,  a)$  is strictly  positive  on  $\R_+$  because of  condition
\reff{eq:conservatif}  and equation  \reff{eq:int_u2}.  Since $\psi $  is
locally  Lipschitz on  $(0,\infty )$,  we deduce  there exists  a unique
locally bounded non-negative solution of \reff{eq:def_w}.
\end{proof}

\subsection{CBI}
\label{sec:CBI}
Let $x>0$,  $\bar \alpha \geq 0$,  $\nu$ be a  Radon measure on  $(0,\infty )$
such  that  $\int_{(0,\infty )}  (1\wedge  x)\;  \nu(dx)<\infty $.   Let
$\cb_+$ denote  the set of  non-negative measurable functions  defined on
$\R$. Let  $h\in \cb_+$ be  locally bounded.  We consider  the following
independent processes.
\begin{itemize}
   \item $\sum_{i\in I} \delta_{t^i,x^i,Z^i}$,  a Poisson point measure with
     intensity $h(t)\ind_{\{t\geq 0\}} dt\; \nu(dx)\; \P_{x,t}(dZ)$.
   \item $\tilde Z$,  distributed according to  $\P_x$. 
    \item $\sum_{j\in J} \delta_{t^j,\hat Z^j}$,  a Poisson point measure with
     intensity $\bar \alpha  h(t)\ind_{\{t\geq 0\}} dt \;Q_t(dZ)$.
\end{itemize}

For $t\in \R$,  let $Y_t=\tilde Z_t+ \sum_{i\in I}  Z_t^i+ \sum_{j\in J}
\hat Z^j_{t}\in  [0,\infty ]$.  We say $Y=(Y_t, t\geq
0)$  is a  continuous  state branching  process  with immigration  (CBI)
started at $x$,  whose branching mechanism is $\psi $  and immigration is
characterized with $(h,\phi)$ where the immigration mechanism, $\phi$,  is defined by
\begin{equation}
   \label{eq:def-phi}
\phi(\lambda) = \bar \alpha  \lambda + \int_{(0,\infty )} \nu(dx) (1-
\expp{-\lambda  x}), \quad \lambda\geq 0 ,
\end{equation}
where $\bar \alpha \geq 0$,  and $\nu$ is a  Radon measure on  $(0,\infty )$
such  that  $\int_{(0,\infty )}  (1\wedge  x)\;  \nu(dx)<\infty $. 

One gets  $Y$ is a  conservative Hunt process  when $h$ is  constant, see
\cite{kw:bpirll}.  Notice that $Y$ is a non-homogeneous Markov
processes. We also have   $Y_0=x$,  and  $Y_t=0$  for $t<  0$.   

Using  Poisson point  measure  property,  one  can  construct   on  the  same
probability space two CBI, $Y^1$  and $Y^2$, with same branching process
$\psi $,  same starting  point  and immigration  characterized by  $(h^1,
\phi)$ and $(h^2,\phi)$ such that $Y^1_t\leq Y^2_t$ for all $t\leq T$ as
soon as $h^1(t)\leq  h^2(t)$ for all $t\leq T$.  We  can apply this with
$h^1=h$ and  $h^2(t)= \sup\{h(s); s\in  [0,T]\}$ for $t\in \R$  and some
$T>0$, and use that $Y^2$  is conservative (see \cite{kw:bpirll}) to get
that  $Y^1$ has a  locally bounded  version over  $[0,T]$. Since  $T$ is
arbitrary, this  implies that any CBI  has a locally  bounded version.  We shall
work with this version.

The following Lemma  is a direct consequence of  the exponential formula
for Poisson point measures (see eg. \cite{ry:mb}, chap. XII).
\begin{lem}
  Let  $\mu$  be a  finite
  positive  measure on  $\R$  with  support bounded from above (i.e. 
  $H^\mu$  is finite).   We have  for
  $s\in \R$:
\begin{equation}
\E\left[\expp{-\int Y_{r-s} \; \mu(dr)
  }\right]=\expp{- xw(s)- \int_0^\infty  h(t) \phi(w(s+t)) dt}, 
\end{equation}
where the function $w$ is defined by \reff{eq:variante-mu}.
\end{lem}


\section{State dependent immigration}
\label{sec:(Y,fi)}
\subsection{Induction formula}
\label{sec:induction}
Let  $(x_k,k\in \N)$  be a  sequence of  non-negative real  numbers. Let
$Y^0$   be  a  CB   with  branching   mechanism  $\psi^0$,   defined  by
\reff{eq:def-psi0},  starting at $x_0$.  We shall  assume that  $Y^0$ is
conservative,   that  is   condition  \reff{eq:conservatif}   holds  for
$\psi^0$.   We construct  by  induction  $Y^n$, $n\geq  1$,  as the  CBI
started  at $x_n$,  with  branching mechanism  $\psi^0$ and  immigration
characterized by $(Y^{n-1},\phi)$, with $\phi$ given by
\reff{eq:def-phi}. 
\begin{lem}
\label{lem:Y0n}
Let $(\mu_k, k\in \N)$ be a family of finite measures on $\R$ with support
bounded from above. 
   We have for all $n\in \N$, $s\in \R$, 
\[
\E\left[\expp{-\sum_{k=0}^n \int Y^k_{r-s}\;  \mu_k(dr)
    }\right]=\expp{-\sum_{k=0}^n x_{n-k}w_k^{(n)} (s)},
\]
where $w_0^{(n)}$  is defined by \reff{eq:variante-mu} with  $\mu$ replaced by
$\mu_n$, and for $k\geq  1$, $w_{k}^{(n)}$ is defined by \reff{eq:variante-mu}
with  $\mu$   replaced  by  $\mu_{n-k}(dr)  +   \phi(w_{k-1}(r))\;  dr$.  In
particular,   $w_k$  is   a  locally bounded non-negative 
solution of the equation
\begin{equation}
   \label{eq:def_w_k}
w(s)+\int_s^{\infty}\psi^0(w(r))dr=\int_{[s,\infty )} \mu_{n-k}(dr)+
\int_s^\infty  \phi(w_{k-1}^{(n)}(r))\; dr, \quad s\in
\R. 
\end{equation}
(Notice $w_k(s)=0$ for $s>\max\{H^{\mu_{k'}}, k'\in \{0, \ldots, k\}\}$.)
\end{lem}
\begin{proof}
  This  is a consequence of the computation of  $\displaystyle
  \E\left[\expp{-\sum_{k=0}^n \int Y^k_{r-s}\; 
      \mu_k(dr) }\Big| Y^0,\ldots, Y^{n-1}\right]$, using   Proposition
  \ref{prop:uniq_x}. This also implies that \reff{eq:def_w_k}
  holds. Then, by induction, one deduces from \reff{eq:def_w_k} that
  $w_k$ is locally bounded. 
\end{proof}

\subsection{Convergence of the  total mass process}
\label{sec:def-X}

We  consider the  sequence  $(Y^n,  n\geq 0)$  defined  in the  previous
Section  with  $x_0=x\geq   0$  and  $x_n=0$  for  $n\geq   1$.  We  set
$X^{n}_t=\sum_{k=0}^n Y^k_t$ for $t\in  \R$. Let $X_t$ be the increasing
limit of $X^n_t$ as $n\to+\infty$, for all $t\in \R$.  We have $X_t\in [0,
+\infty ]$.  We call $X=(X_t, t\in  \R)$ a CBI  with branching mechanism
$\psi^0$ and immigration process $(X,\phi)$. We set $\psi=\psi^0-\phi$. 

\begin{rem}
\label{rem:phi}
For $\lambda\geq 0$, we have 
\[
\psi^0(\lambda) -\phi(\lambda)=\Big(\alpha^0 -\bar \alpha - \int_{(0,1]}\ell\;
\nu(d\ell)
 \Big)\lambda + \beta \lambda^2
+ \int_{(0,\infty )} (\pi(d\ell)+\nu(d\ell)) 
\left[\expp{-\lambda\ell}-1+\lambda\ell\ind_{\{\ell \leq 1\}}\right].  
\]
This gives  that $\psi=\psi^0 -\phi$ is  a branching mechanism.
\end{rem}

The  process $Y^0$  describes the  size process  of  the Eve-population,
$Y^1$  the size  process  of the  population  of mutants  born from  the
Eve-population  $Y^0$,  $Y^2$ the  size  process  of  the population  of
mutants born from  mutant population $Y^1$, and so  on. The size process
of the total  population is given by $X=\sum_{k\geq  0} Y^k$. In neutral
mutation case, it is natural to assume that all the processes $Y^k$ have
the same branching mechanism. Since we assume $x_k=0$ for all $k\geq 1$,
this means only the Eve-population is present at time $0$.

\begin{theo}
\label{th:main}
We    assume    that   $\psi$    is    conservative,   i.e.    satisfies
\reff{eq:conservatif}.  The  process $X$, which is a  CBI with branching
mechanism  $\psi^0$ and  immigration process  $(X,\phi)$, is  a  CB with
branching mechanism $\psi=\psi^0-\phi$.
\end{theo}

\begin{rem}
   \label{rem:X=M}
   As a consequence of Theorem \ref{th:main}, $X$ is a Markov process.
   Notice that $(Y^0,  \ldots, Y^n)$ is also Markov  but not
   $(X^n_t,t\ge 0)$ for
   $n\geq 1$.
 \end{rem}

\begin{proof}
Let $\mu$ be a finite measure on $\R$ with support bounded
from above (i.e. $H^\mu<\infty )$. We shall assume that
$\mu(\{H^\mu\})=a>0$. 

We keep the notations of Lemma \ref{lem:Y0n}, with $\mu_k=\mu$. In
particular we see from \reff{eq:def_w_k} that $w_k^{(n)} $ does not
depend on $n$. We shall denote it by $w_k$. 
By monotone convergence, we have 
\[
\E\left[\expp{- \int X_{r-s}\;  \mu(dr)
    }\right]=\lim_{n\rightarrow \infty } \E\left[\expp{- \sum_{k=0}^n
      \int Y^k _{r-s}\;  \mu(dr)
    }\right]= \lim_{n\rightarrow \infty }\expp{-xw_n(s)},
\]
where the limits are non-increasing.   This implies that $(w_n, n\geq 0)$
increases  to   a  non-negative   function  $w_\infty  $.   By  monotone
convergence  theorem (for $\int_s^{H^\mu}\psi^0(w(r))\ind_{\{w_n(r)>0\}}\;
dr$ and the integral with $\phi$) and dominated convergence theorem (for
$\int_s^{H^\mu}\psi^0(w(r))\ind_{\{w_n(r)\leq 0\}}\;  dr$), we deduce from
\reff{eq:def_w_k}, that $ w_\infty $ solves $ w(s)=0$ for $s> H^\mu$ and
\begin{equation}
   \label{eq:bar_w}
w(s)+\int_s^{H^\mu}\psi^0(w(r))dr=\int_{[s,\infty )} \mu(dr)+
\int_s^{H^\mu} \phi(w(r))\; dr, \quad s\leq H^\mu.
\end{equation}
Notice that $w_\infty (s)\in [0,\infty ]$ and the two sides of the previous
equality may be infinite. 

Thanks  to Proposition  \ref{prop:uniq_x},  and since  $\psi^0-\phi$ is  a
branching mechanism  (see Remark  \ref{rem:phi}), there exists  a unique
locally bounded non-negative solution of \reff{eq:bar_w}, which we shall
call $\bar w$. Therefore to prove  that $w_\infty =\bar w$, it is enough
to check that  $w_\infty $ is locally bounded. This will  be the case if
we check  that $w_\infty \leq \bar  w$. In particular,  we get $w_\infty
=\bar w$, if we  can prove that $w_n\leq \bar w $  for all $n\in \N$. We
shall prove this by induction.

We consider the measure $\mu^0(dr)=\mu(dr)+ \phi(\bar w(r))\ind_{\{r\leq
  H^\mu\}}\; dr$. Notice $H^{\mu_0}=H^\mu$ and
$\mu^0(\{H^{\mu^0}\})=\mu(\{H^{\mu}\})=a>0$. 
We define $\bar w_0$ by 
\[
\expp{-x \bar w_0(s)}=\E\left[\expp{- \int Y^0_{r-s}\;  \mu^0(dr)
    }\right].
\]
The function  $\bar w_0$ is a locally bounded non-negative function which solves 
\[
w(s)+\int_s^{H^\mu}\psi^0(w(r))dr=\int_{[s,\infty )} \mu(dr)+
\int_s^{H^\mu} \phi(\bar w(r))\; dr, \quad s\leq H^\mu.
\]
Thanks to Proposition \ref{prop:uniq_x}, $\bar w_0$ is unique. Since $\bar w$
solves the same equation, we deduce that $\bar w_0=\bar w$.  
We also have 
\[
\expp{-x w_0(s)}=\E\left[\expp{- \int Y^0_{r-s}\;  \mu(dr)
    }\right]
\geq  \E\left[\expp{- \int Y^0_{r-s}\;  \mu^0(dr)
    }\right]. 
\]
This implies that $w_0\leq \bar w_0=\bar w$. 

Assume we proved that $w_{n-1}\leq \bar w$ for some $n\geq 1$. Then we can
consider the measure $\mu^n(dr)=\mu(dr)+ [\phi(\bar w(r))-
\phi(w_{n-1}(r))]\ind_{\{r\leq 
  H^\mu\}}\; dr$. Notice $H^{\mu^n}=H^\mu$ and
$\mu^n(\{H^{\mu^n}\})=a>0$. 
Recall $x=x_0\geq 0$ and $x_k=0$ for $k\geq 1$. We define $\bar w_{n}$ by 
\[
\expp{-x \bar w_n(s)}=\E\left[\expp{- \int Y^{n}_{r-s}\;  \mu^{n}(dr)
    }\right].
\]
The function  $\bar w_n$ is a locally bounded non-negative function which
solves  for $ s\leq H^\mu $
\begin{align*}
w(s)+\int_s^{H^\mu}\psi^0(w(r))dr
&=\int_{[s,\infty )} \mu^n(dr)+
\int_s^{H^\mu} \phi(w_{n-1} (r))\; dr\\
&=\int_{[s,\infty )} \mu(dr)+
\int_s^{H^\mu} \phi(\bar w (r))\; dr.
\end{align*}

Thanks to Proposition \ref{prop:uniq_x}, $\bar w_n$ is unique. Since $\bar w$
solves the same equation, we deduce that $\bar w_n=\bar w$.  
We also have 
\[
\expp{-x w_n(s)}=\E\left[\expp{- \int Y^n_{r-s}\;  \mu(dr)
    }\right]
\geq  \E\left[\expp{- \int Y^n_{r-s}\;  \mu^n(dr)
    }\right]. 
\]
This implies that $w_n\leq \bar w$. Therefore, this holds for all $n\geq
0$, which according to our previous remark entails that $w_\infty =\bar
w$.

By taking $\mu(dr)=\sum_{k=1}^K \lambda_k \delta_{t_k}(dr)$ for $K\in
\N^*$, $\lambda_1, \ldots, \lambda_K\in [0,\infty )$ and $0\leq t_1\leq
\ldots \leq  t_K$, we deduce that $X$ has the same finite marginals
distribution   as a CB with branching mechanism $\psi^0-\phi$. Hence $X$
is a CB with branching mechanism $\psi^0-\phi$.
   
\end{proof}

\section{The dual to the pruning at node}
\label{sec:dual}
For $\theta\in \R$, we consider the group of operators  $(T_\theta,
\theta\in \R)$  on the set of
real measurable functions defined by
\[
T_\theta(f)(\cdot)=f(\theta+\cdot) -f(\theta).
\]
Let  $\psi^0$ be given  by \reff{eq:def-psi0}  with Lévy  measure $\pi$.
Using the previous Section, for  $\theta>0$, we can give a probabilistic
interpretation to  $T_{-\theta}(\psi^0)$ as  a branching mechanism  of a
CBI  with  proportional  immigration.  Let  $\theta_0=  \sup\{\theta\geq
0;\quad\int_{(1,\infty )}  \expp{\theta \ell} \;  \pi(d\ell)<\infty\} $.
Notice    that   $\theta_0=0$    if   ${\psi^0}'(0^+)=-\infty    $,   as
${\psi^0}'(0^+)=-\infty  $ is equivalent  to $  \int_{(1,\infty )}\ell\;
\pi(d\ell)=+\infty   $.     We   assume   $\theta_0>0$    and   we   set
$\Theta=(0,\theta_0]$  if $\int_{(1,\infty  )}  \expp{\theta_0 \ell}  \;
\pi(d\ell)<\infty$ and  $\Theta=(0,\theta_0)$ otherwise.  Let $\theta\in
\Theta$.  We define
\[
\phi_\theta(\lambda) = 2\beta \theta\lambda+ \int_{(0,\infty )} (\expp{\theta
  x}-1)(1- \expp{-\lambda  x})  \; \pi(dx) .
\]
It     is     straightforward     to     check     that     $T_{-\theta}
(\psi^0)=\psi^0-\phi_\theta$  and that  $\phi_\theta$ is  an immigration
mechanism.  Notice  that  for  $\theta<\theta_0$, we  have  $T_{-\theta}
(\psi^0)'(0^+)>-\infty   $,  that   is  $T_{-\theta}   (\psi^0)$   is  a
conservative branching mechanism.

The next Corollary is a direct
consequence of the previous Section. 
\begin{cor}
\label{cor:CBI-CB}
  Let $\theta\in \Theta$. If $\theta=\theta_0$ assume furthermore that
  $T_{-{\theta_0}}({\psi^0})$ is conservative.  
A CBI process  $X$  with branching mechanism ${\psi}$ and immigration
  $(X, \phi_\theta )$ is a CB with branching mechanism
  $T_{-\theta}({\psi})$. 
\end{cor}

On the  other end, for $\theta>0$,  $T_\theta(\psi)$ can be  seen as the
branching mechanism of a  pruned CB. The following informal presentation
relies on the pruning procedure developed in \cite{ad:falp}.
Let  us  consider a  CRT  associated  with  a critical  or  sub-critical
branching mechanism $\psi$  with no Brownian part, which  we shall write
in the following form:
\[
\psi(\lambda)=\alpha_1\lambda+\int_{(0,+\infty)}(e^{-\lambda
  r}-1-\lambda r)\pi(dr),
\]
with     $    \int_{[1,\infty    )}     r\;    \pi(dr)     <\infty    $,
$\alpha_1=\alpha+\int_{[1,\infty   )}    r\;   \pi(dr)   \ge    0$   and
$\int_{(0,1)}r\pi(dr)=+\infty$.  In  that  case,  the  CB  process  with
branching  mechanism  $\psi$ has   no  diffusion  part  ($\beta=0$)  and
increases only by positive jumps.
Let us  recall that a CRT can be coded  by the so called
height process  $H=(H_t,t\geq 0)$, see  \cite{dlg:rtlpsbp}. Intuitively,
for the individual $t\geq 0$, $H_t\geq 0$ represents its generation. The
individual  $t$ is  called an  ancestor of  $s$ if  $H_t=\min\{H_u, u\in
[s\wedge  t,  s  \vee  t]\}$,  and we  shall  write  $s\succcurlyeq  t$.
Informally, for $t$  fixed, the ``size'' of the  population at generation
$a\geq 0$ of all individuals $r\leq t$ is given by the local time of $H$
at level  $a$ up to time $t$,  $L^a_t$ say. For the  CRT associated with
the branching  mechanism $\psi$,  the process $L=(L^a_{T_x},  a\geq 0)$,
where  $T_x=\inf\{t\geq  0,  L^0_t\geq  x\}$  is  a  CB  with  branching
mechanism $\psi$. The  height process codes for the  genealogy of the CB
process $L$.

An individual  $t$ is  called a node  of the  CRT if the  height process
corresponding  to its  descendants, $(H_s-H_t,  s\succcurlyeq t)$  has a
positive local time at level  $0$, say $\Delta_t$. (If $t\leq T_x$, then
$\Delta_t$ corresponds  to a jump of  the CB process $L$  at level $H_t$;
reciprocally to a jump $\Delta$ of the CB process $L$ at level $a$ there
corresponds   an  individual   $t\leq   T_x$  such   that  $H_t=a$   and
$\Delta_t=\Delta$.)  Intuitively $\Delta_t$  corresponds to the ``size''
of the offspring population of individual $t$.  Let $\theta>0$ be fixed.
A   node   $t$  of   size   $\Delta_t$   is   marked  with   probability
$1-\expp{-\theta\Delta_t}$, independently of  the other nodes.  To prune
the CRT, we just remove all  individuals who have a marked ancestor. The
height process of  the pruned CRT is then given  by $H^\theta=(H_{C_t}, t\geq
0)$  where $C$ is  the inverse  of the  Lebesgue measure  of the  set of
individuals whose ancestors have no mark:
\[
 C_t=\inf  \{r_0\geq 0;  \int_0^{r_0} \ind_{\{\forall s,
r\succcurlyeq s, 
s\text{  is not  marked}\}}\; dr  \geq t\}.
\]
Theorem 6.1 in  \cite{ad:falp} shows  that this
pruned CRT is itself a CRT associated with the branching mechanism
$T_\theta(\psi)$.

By  looking at  the local  time of  the pruned  process, we  get  a nice
construction of a CB process of branching mechanism $T_\theta (\psi)$, which
we shall called a pruned CB with intensity $\theta>0$, from a CB process
of branching mechanism $\psi$.   Notice this construction was done under
the assumption that $\beta=0$ (see also \cite{as:psf} when $\beta>0$ and
$\pi=0$).  The  general pruning procedure in the  general case $\beta>0$
and   $  \pi\neq   0$  will   be  presented   in  a   forthcoming  paper
\cite{adv:plcrt}.

In a certain sense the immigration is the dual to the pruning at node:
to  build   a  CB  process of  branching
mechanism $\psi$ from  a CB process of branching mechanism
$T_\theta(\psi)$, with $\theta>0$, 
one has to  add  an immigration at  time $t$ which rate is proportional to the
size of the population at time~$t$ and immigration mechanism $\tilde
\phi_\theta$ defined by: 
\[
\tilde \phi_\theta(\lambda)= T_\theta(\psi)(\lambda) -\psi(\lambda)= 
2\beta \theta\lambda+ \int_{(0,\infty )} (1-\expp{-\theta
  x})(1- \expp{-\lambda  x})  \; \pi(dx), \quad \text{for $\lambda\geq 0$}.
\]
In other words, we get the  following result, whose first part comes from
Theorem 6.1 in \cite{ad:falp}. As  in \cite{ad:falp}, we assume only for
the   next  Corollary   that   $\beta=0$  and   $\int_{(0,1)}  \ell   \;
\pi(d\ell)=+\infty$.

\begin{cor}
\label{cor:dual}
Let  $X$  be  a  critical  or sub-critical  CB  process  with  branching
mechanism  $\psi$.  Let  $X^{(\theta)}$ be  the  pruned CB  of $X$  with
intensity $\theta>0$  : $X^{(\theta)} $  is a CB process  with branching
mechanism $T_\theta(\psi)$.  The CBI process, $\tilde X$, with branching
mechanism   $T_\theta(\psi)$   and   immigration  $(\tilde   X,   \tilde
\phi_\theta)$ is distributed as $X$.
\end{cor}

\section{Application : law of the initial process}
\label{sec:appli}

We consider a population whose size  evolves as $X=(X_t, t\geq 0)$, a CB
with  branching  mechanism  $\psi$ given by \reff{eq:def_psi}.   We
assume  $\psi$  satisfies  the 
hypothesis of  Section \ref{sec:CB-CBI}. This  population undergoes some
irreversible mutations with constant  rate. Each mutation produces a new
type  of individuals.  In  the critical  or sub-critical  quadratic case
($\pi=0$) this corresponds to the  limit of the Wright-Fisher model, but
for the fact that the ``size''  of the population is not constant but is
distributed  as a  CB. 

We assume the population at time $0$ has the same original Eve-type.  We
are interested in the law of $Y^0=(Y^0_t, t\geq 0)$, the ``size'' of the
sub-population  with the  original type  knowing the  size of  the whole
population.  In  particular, we shall compute  $\P(Y^0_t=0| X_t>0)$, the
probability for the Eve-type  to have disappeared, conditionally on the
survival of the  total population at time $t$.

We   shall   assume   $Y^0$   is   a   CB   with   branching   mechanism
$\psi^0$  and  $X$  is  the  CBI  with
immigration $(X,\phi)$,  with $\phi=\psi^0-\psi$, considered  in Section
\ref{sec:def-X}. Thus, we model  the mutations by an immigration process
with   rate  proportional  to   the  size   of  the   population.

The joint law of $(X_t,Y_t^0)$ can be easily characterized by the
following Lemma. 
\begin{lem} 
\label{lem:loi_Xt_Yt}
Let $t\geq 0$, $\lambda_1, \lambda_2 \in \R_+$. We assume
  $X_0=Y^0_0=x\geq 0$. We have 
\[
\E\left[ \expp{-\lambda_1X_t-\lambda_2Y_t^0}\right]=\expp{-x w(0)},
\]
where  $(w,w^*)$  is  the  unique  measurable  non-negative  solution  on
$(-\infty ,t]$ of
\begin{align*}
w(s)+\int_s^t\psi^0\bigl(w(r)\bigr)dr 
&=\lambda_1+\lambda_2+\int_s^t\phi\bigl(w^*(r)\bigr)dr, \\
w^*(s)+\int_s^t\psi\bigl(w^*(r)\bigr)dr 
& =\lambda_1.
\end{align*}
\end{lem}

\begin{proof}
  Recall  notation  of  Section  \ref{sec:def-X}. In particular $x_0=x$
  and $x_n=0$ for all $n\geq 1$.  Let  us  apply  Lemma
  \ref{lem:Y0n} with
\begin{align*}
\mu_0(dr) & =(\lambda_1+\lambda_2)\delta_t(dr),\\
\mu_k(dr) & =\lambda_1\delta_t(dr)\quad\mbox{for }k\ge 1.
\end{align*}
We get
\[
\E\left[\expp{-(\lambda_1X_t^n+\lambda_2Y_t^0)}\right]=\expp{-xw_n^{(n)}(0)},
\]
where  for $s\le t$,
\begin{align*}
w_0^{(n)}(s)+\int_s^t\psi^0\bigl(w_0^{(n)}(r)\bigr)dr & =\lambda_1,\\
w_k^{(n)}(s)+\int_s^t\psi^0\bigl(w_k^{(n)}(r)\bigr)dr &
=\lambda_1+\int_s^t\phi\bigl(w_{k-1}^{(n)}(r)\bigr)dr\quad\mbox{for
}1\le k\le n-1,\\ 
w_n^{(n)}(s)+\int_s^t\psi^0 \bigl(w_n^{(n)}(r)\bigr)dr &
=\lambda_1+\lambda_2+\int_s^t\phi\bigl(w_{n-1}^{(n)}(r)\bigr)dr. 
\end{align*}
We let  $n$ goes to infinity and use similar arguments as in the proof
of Theorem \ref{th:main}  to get the result.
\end{proof}

Some more  explicit computations  can be made  in the case  of quadratic
branching mechanism  (see also \cite{w:bprtsbm}  when $\alpha =0$). Let
$\alpha \geq 0$ and $\theta>0$ and set
\[
\psi(u)=\alpha   u+  u^2,   \quad  \psi^0(u)=T_\theta   (\psi)(u)  =(\alpha
+2\theta) u+ u^2.
\]
The CB which  models the total population is  critical $(\alpha =0$) or
sub-critical    ($\alpha >0$).     The    immigration   mechanism    is
$\phi(u)=\psi^0(u) -\psi(u)=2\theta u$.

We set $b=(\alpha +2\theta)$ and for $t\geq 0$,
\begin{equation}\label{eq:sol_u*}
h(t)=\begin{cases}
1+ \lambda_1  
  \frac{1-\expp{-\alpha 
      t}}{\alpha }& \text{if 
$\alpha >0$,}\\
1+\lambda_1 t & \text{if
$\alpha =0$.}
\end{cases}
\end{equation}

\begin{prop}
\label{prop:loi_Xt_Yt}
Let $t\geq 0$, $\lambda_1, \lambda_2 \in \R_+$. We have
\[
\E\left[
  \expp{-\lambda_1X_t-\lambda_2Y_t^0}\right]
=\expp{-x v_0(t)},
\]
where 
\[
v_0(t)= \expp{-bt} h(t)^{-2} \left(\frac{1}{\lambda_2}+\int_0^{t}\expp{-br}
  h(r)^{-2}  dr\right)^{-1} 
+ \lambda_1 \expp{-\alpha  t}  h(t)^{-1}  .
\]
\end{prop}

\begin{proof}
By the previous lemma, we have
\begin{equation}
   \label{eq:v0}
\E\left[ \expp{-\lambda_1X_t-\lambda_2Y_t^0}\right]=\expp{-xw(0)}
\end{equation}
where for $s\le t$,
\begin{align}
\label{eq:w}
&w(s)+\int_s^tw(r)\bigl(w(r)+b\bigr)dr 
=\lambda_1+\lambda_2+2\theta \int_s^t w^*(r)dr,\\
\nonumber &w^*(s)+\int_s^tw^*(r)\bigl(w^*(r)+\alpha \bigr)dr  =\lambda_1.
\end{align}
The last equation is equivalent to
\begin{equation}
\label{eq:w*}
(w^*)'-w^*(w^*+\alpha )=0 \quad\text{on $(-\infty , t]$}, \quad
w^*(t)=\lambda_1.
\end{equation}
The function $\displaystyle z^*:=\frac{1}{w^*}$ is thus the unique solution of 
\[
(z^*)'+\alpha  z^*+1=0\quad\text{on $(-\infty , t]$},\quad 
z^*(t)=\frac{1}{\lambda_1}.
\]
If $\alpha >0$, this  leads to
\[
z^*(s)
=\frac{1}{\alpha }\left(\expp{\alpha (t-s)}-1\right)
+\frac{1}{\lambda_1}\expp{\alpha  (t-s)}. 
\]
If $\alpha =0$, we have $\displaystyle z^*(s)=
t-s+\inv{\lambda_1} $. We get
\begin{equation}\label{eq:sol_w*}
   w^*(s)=h'(t-s)h(t-s)^{-1}=\lambda_1 \expp{-\alpha  (t-s)}  h(t-s)^{-1} .
\end{equation}
Equation \reff{eq:w} is equivalent to
\[
w'-w(w+b)=-2\theta w^*\quad\text{on $(-\infty , t]$}, \quad 
w(t)=\lambda_1+\lambda_2.
\]
Set $y=w-w^*$ and use the differential equation \reff{eq:w*}, to get 
that $y$ solves 
\[
y'-y^2 -y(2w^*+b)=0\quad\text{on $(-\infty , t]$}, \quad y(t)=\lambda_2.
\]
Then the function $z:=1/y$ is the unique solution of
\[
z'+(2w^*+b)z+1=0\quad\text{on $(-\infty , t]$}, \quad 
z(t)=\frac{1}{\lambda_2}.
\]
One solution of the homogeneous differential equation
$z_0'=-(2w^*+b)z_0$
is $z_0(s)=\expp{b(t-s)}h(t-s)^2$. 
Looking for solutions of the form $z(s)=C(s)z_0(s)$ gives
$$z(s)=z_0(s)\left(\frac{1}{\lambda_2}+\int_s^{t}z_0(u)^{-1} du\right).$$
We conclude using \reff{eq:v0} and $w=w^*+z^{-1}$. 
\end{proof}

\begin{rem}
   \label{rem:Yt-cond}
We can compute the conditional probability of the non extinction of
the Eve-population: $\P(Y^0_t>0|X_t>0)$. However, this computation can
be done without the joint law of $(X_t, Y^0_t)$ as 
\[
\P(Y^0_t>0|X_t>0)=
\frac{\P(Y^0_t>0,X_t>0)}{\P(X_t>0)}=\frac{\P(Y^0_t>0)}{\P(X_t>0)}
=\frac{1-\P(Y^0_t=0)}{1-\P(X_t=0)} , 
\]
with  $\displaystyle \P(X_t=0)=\lim_{\lambda_1\rightarrow \infty }
\E[\expp{-\lambda_1 
  X_t}]= \expp{-xg(\alpha ,t)^{-1}}$ and $\displaystyle 
\P(Y^0_t=0)=\lim_{\lambda_2\rightarrow 
  \infty } \E[\expp{-\lambda_2 
Y^0_  t}]=\expp{-x g(b,t)^{-1}}$, where 
\begin{equation}
   \label{eq:g}
g(a,t)=\begin{cases}
  \frac{\expp{a
      t} -1}{a}& \text{if 
$a>0$,}\\
 t & \text{if
$a=0$.}
\end{cases}
\end{equation}
\end{rem}
The same kind of computation allows also to compute the joint law at
different times.
\begin{prop}
\label{prop:v1}
Let $0\leq u<t$, $\lambda_1,\lambda_2 \in \R_+$. We have
\[
\E\left[\expp{-\lambda_1X_t-\lambda_2Y_u^0}\right]=\expp{-x v_1(u,t)},
\]
where  
\[
v_1(u,t)= \expp{-bt} h(t)^{-2} \left(\frac{\expp{-b(t-u)} h(t-u)^{-2}  }{\lambda_2}+\int^t_{t-u}\expp{-br}
  h(r)^{-2}  dr\right)^{-1} 
+ \lambda_1 \expp{-\alpha  t}  h(t)^{-1}  .
\]
\end{prop}

\begin{proof}
 Recall  notation  of  Section  \ref{sec:def-X}. In particular $x_0=x$
  and $x_n=0$ for all $n\geq 1$.  Let  us  apply  Lemma
  \ref{lem:Y0n} with
\begin{align*}
\mu_0 & =\lambda_1\delta_t+\lambda_2\delta_u,\\
\mu_k & =\lambda_1\delta_t\qquad\mbox{for }k\ge 1.
\end{align*}
Let $n$ goes to infinity as in the proof of Lemma \ref{lem:loi_Xt_Yt} to
get that 
\begin{equation}
   \label{eq:v1}
\E\left[\expp{-\lambda_1X_t-\lambda_2Y_u^0}\right]=\expp{-xw(0)},
\end{equation}
where  $(w,w^*)$ is the unique non-negative solution on $(-\infty ,t]$  of
\begin{align}
\label{eq:wbis}w(s)+\int_s^t\psi^0\bigl(w(r)\bigr)dr &
=\lambda_1+\lambda_2\ind_{\{s\le u\}}+\int_s^t\phi\bigl(w^*(r)\bigr)dr,\\ 
\nonumber
w^*(s)+\int_s^t\psi\bigl(w^*(r)\bigr)dr &
=\lambda_1.
\end{align}
Notice $w^*$ is still given by  \reff{eq:sol_w*}. 
For $s>u$, we have $w(s)=w^*(s)$ and, for $s\le u$, Equation
\reff{eq:wbis} is equivalent to 
\[
w'-w(w+b)=-2\theta w^*\quad\text{on $(-\infty , u]$}, \quad 
w(t)=w^*(u)+\lambda_2.
\]
{From} the proof of Proposition \ref{prop:loi_Xt_Yt}, we get 
\[
\frac{1}{w(s)-w^*(s)}
=\expp{b(t-s)}h(t-s)^2 
\left(\frac{\expp{-b(t-u)}h(t-u)^{-2}}{\lambda_2}
  +\int_s^u\expp{-b(t-r)} h(t-r)^{-2}
dr\right).   
\]
We conclude using \reff{eq:v1}. 
\end{proof}

At this stage, we can give the joint distribution of the extinction time
of $X$,  $\tau_X=\inf \{t>0; X_t=0\}$, and  of $Y^0$, $\tau_{Y^0}=\inf\{
t>0;  Y^0_t=0\}$.   For $u\leq  t$,  we  have $\P(\tau_X\leq  t,\tau_{Y^0}\leq
u)=\lim_{\lambda_1 \rightarrow \infty  , \; \lambda_2 \rightarrow \infty
} \exp - x v_1(u,t)$ that is
\[
\P(\tau_X\leq t,\tau_{Y^0}\leq u)= \exp - x \left(   \expp{-bt+2\alpha  t }
  (\int _{t-u}^t \!\!\!\!\expp{-br+2\alpha  r} g(\alpha ,t)^2
  g(\alpha ,r)^{-2} 
  dr)^{-1} + g(\alpha ,t)^{-1} \right). 
\]
We can compute the probability of simultaneous extinction of the
Eve-population and the whole population, see also proposition 5 in
\cite{w:bprtsbm}, where $\alpha =0$.  In \cite{ad:wdlcrtseppnm}, using
different techniques we derive this formula for general critical or
sub-critical branching mechanisms. 

\begin{lem}
   \label{lem:tau=sigma}
We have $\displaystyle \P(\tau_{Y^0}=\tau_X|\tau_X=t)=\expp{-2\theta t}$.
\end{lem}

\begin{proof}
   We have $\displaystyle 
\P(\tau_{Y^0}=\tau_X|\tau_X=t)=1-\frac{\lim_{u \uparrow t} \partial_t
  \P(\tau_{Y^0}\leq u, \tau_X\leq t)}{\partial_t
  \P(\tau_X\leq t)}=\expp{-2\theta t}$.
\end{proof}

We can deduce from the latter Proposition the law of $Y_u^0$
conditionally on the  non-extinction of the whole population. 
We set 
\[
A(b,u)=\frac{1}{\lambda_2}\expp{bu}+g(b,u).
\]

\begin{prop}
\label{prop:Yu_Xi}
Let $u\geq 0$, $\lambda_2\in \R_+$. We have 
\[
\lim_{t\to+\infty}\E\left[\expp{-\lambda_2Y_u^0}\Bigm|X_t>0\right]=\expp{-xA(b,u)^{-1} 
  } \left(1-A(b,u)^{-2}G(\alpha ,u)
\right),
\]
where   
\[
G(a,u)=\frac{2}{\lambda_2} \expp{bu}
  g(a,u)+ \begin{cases}
  2\frac{g(b+a,u)- 
  g(b,u)}{a}  & \text{if 
$a>0$,}\\
2\partial_1
  g(b,u)  & \text{if
$a=0$.}
\end{cases}
\]
\end{prop}

\begin{proof}
We have 
\[
\E\left[\expp{-\lambda_2Y_u^0}\biggm|X_t>0\right]=\frac{\E\left[\expp{-\lambda_2Y_u^0}\right]-\E\left[\expp{-\lambda_2Y_u^0}\ind_{\{X_t=0\}}\right]}{\P(X_t>0)}\cdot
\]
Using Proposition \ref{prop:v1}
\[
\E\left[\expp{-\lambda_2Y_u^0}\ind_{\{X_t=0\}}\right]  =
\lim_{\lambda_1\to+\infty}\E\left[\expp{-\lambda_2Y_u^0-\lambda_1X_t}\right]=
\expp{ - x \bar v_1(u,t)},
\]
with $\bar v_1(u,t)=\lim_{\lambda_1\to+\infty} v_1(u,t)$.

Definition \reff{eq:g} implies
\[
\bar  v_1(u,t)= 
  \left(\frac{\expp{bu} g(\alpha ,t)^2}
{\lambda_2g(\alpha ,t-u)^{2}}+\expp{bt} \int^t_{t-u}\expp{-br}
  \frac{g(\alpha ,t)^2}{g(\alpha ,r)^{2}}  dr\right)^{-1} 
+  \expp{-\alpha  t}  g(\alpha ,t)^{-1}  .
\]
Performing an asymptotic expansion of $\bar v_1$ as $t$ goes to
$\infty$ leads to the result.
\end{proof}

\newcommand{\sortnoop}[1]{}

\end{document}